\documentclass[11pt,reqno]{amsart}
\usepackage{amsmath,amssymb,amsfonts,amsthm}
\IfFileExists{mathtools.sty}{\usepackage{mathtools}}{}
\usepackage[margin=1.22in]{geometry}
\IfFileExists{enumitem.sty}{%
  \usepackage{enumitem}
  \setlist[itemize]{leftmargin=2em,itemsep=0.25em,topsep=0.35em}
  \setlist[enumerate]{leftmargin=2.2em,itemsep=0.25em,topsep=0.35em}
}{}
\IfFileExists{xcolor.sty}{\usepackage{xcolor}}{\usepackage{color}}
\IfFileExists{microtype.sty}{\usepackage{microtype}}{}
\usepackage[colorlinks=true,linkcolor=blue,citecolor=blue,
            urlcolor=blue]{hyperref}
\numberwithin{equation}{section}
\allowdisplaybreaks

\newtheorem{thm}{Theorem}[section]
\newtheorem{prop}[thm]{Proposition}
\newtheorem{lemma}[thm]{Lemma}

\theoremstyle{definition}

\theoremstyle{remark}

\newcommand{\tr}{\operatorname{tr}}

\theoremstyle{plain}
\newtheorem*{claim}{Claim}
\newcommand{\bbR}{\mathbb R}
\newcommand{\bbS}{\mathbb S}
\newcommand{\cC}{\mathcal C}
\newcommand{\cK}{\mathcal K}
\newcommand{\cL}{\mathcal L}

\title[Capillary $L_p$ Dual Minkowski Flow]{Capillary $L_p$ Dual Minkowski Flow}
\date{}
\author[J. Hu]{Jinrong Hu}
\address{Institut f\"{u}r Diskrete Mathematik und Geometrie, Technische Universit\"{a}t Wien, Wiedner Hauptstrasse 8-10, 1040 Wien, Austria
 }
\email{jinrong.hu@tuwien.ac.at}

\author[Y. Hu]{Yingxiang Hu}
\address{School of Mathematical Sciences, Beihang University, Beijing 100191, China}
\email{huyingxiang@buaa.edu.cn}

\author[B. Yang]{Bo Yang}
\address{Department of Mathematical Sciences, Tsinghua University, Beijing 100084, P.R. China}
\email{ybo@tsinghua.edu.cn}
\address{Institut f\"{u}r Diskrete Mathematik und Geometrie, Technische Universit\"{a}t Wien, Wiedner Hauptstra{\ss}e 8-10, 1040 Wien, Austria}
\email{bo.yang@tuwien.ac.at}
\begin{document}

\begin{abstract}
For $n\geq2$ and $\theta\in(0,\pi/2)$, we study the long-time existence and asymptotic behavior of an even normalized capillary $L_p$ dual Minkowski flow in the Euclidean half-space $\overline{\mathbb{R}^{n+1}_+}$. When $pq\geq0$, the flow exists for all positive time and the evolving hypersurfaces converge smoothly along a sequence of times tending to infinity to a stationary solution of the normalized problem. If $p\neq q$, a dilation of the limit solves the unnormalized problem.
\end{abstract}


\thanks{J. Hu was supported by the Austrian Science Fund (FWF): 10.55776/ESP1358925. Y. Hu was supported by the National Key Research and Development Program of China 2021YFA1001800 and NSFC grant 12671060. B. Yang was supported by China Postdoctoral Science Foundation No.~2024M751605.} \maketitle

\section{Introduction}\label{Sec1}
The Minkowski problem and its extensions connect the geometry of convex bodies with fully nonlinear curvature equations. In the capillary setting, the hypersurface has a boundary constrained to a supporting hyperplane and meets that hyperplane at a prescribed contact angle. Prescribing curvature data therefore leads to a boundary value problem in which the contact-angle constraint must be incorporated into the analysis. Mei, Wang and Weng \cite{MWW25a} introduced and solved the capillary Minkowski problem in the smooth setting, establishing a capillary counterpart of the classical prescribed Gauss curvature problem. Their formulation leads to a Monge-Amp\`ere equation on a spherical cap with a Robin boundary condition.

This formulation has led to several extensions of the capillary Minkowski problem. For acute contact angles, Mei, Wang and Weng \cite[Theorem~1.1]{MWW25b} proved existence and uniqueness for general positive smooth data when $p>n+1$, a normalized eigenvalue result with uniqueness up to dilation when $p=n+1$, and smooth even solutions for even data when $1<p<n+1$. Hu and Ivaki \cite[Theorem~1.2]{HI25a} solved the even problem for $-n-1<p<1$ by iterating capillary curvature images. Here the ambient dimension is $n+1$. Related developments include work on the capillary Orlicz-Minkowski problem \cite{WZ25,LL26} and on the capillary Christoffel-Minkowski problems prescribing area measures \cite{HIS25,MWW26}. These developments extend several Minkowski-type problems in Euclidean space (see, e.g.,  \cite{M897, M903, A38, A39, N53, CY76, P78,L93,HLYZ10}) to the capillary setting, providing a framework for studying Minkowski-type problems with a capillary boundary, while also highlighting the role of boundary estimates.

The dual Minkowski problem, introduced by Huang, Lutwak, Yang, and Zhang \cite{HLYZ16}, concerns the prescription of dual curvature measures of convex bodies, while the $L_p$ dual curvature measures were subsequently introduced by Lutwak, Yang, and Zhang \cite{LYZ18}. These developments incorporate radial information into the framework for prescribing curvature measures, leading to the $L_p$ dual Minkowski problem, which encompasses both the  $L_p$ Minkowski problem and the dual Minkowski problem as special cases. In the setting of closed hypersurfaces, besides the variational approaches developed in \cite{HLYZ16, HZ18}, geometric flows also offer an alternative and effective method for studying dual Minkowski-type problems. See, for example,  Chen-Huang-Zhao \cite{CHZ19}, Li-Sheng-Wang \cite{LSW20} and Chen-Li \cite{CL21}. In particular, the authors of \cite{CHZ19} used anisotropic Gauss curvature flows to obtain smooth even solutions to the $L_p$ dual Minkowski problem for $pq\geq0$, with a normalized formulation in the case $p=q$. This provides a natural motivation for investigating the corresponding  prescribed curvature problem for capillary hypersurfaces.

Given a smooth positive function $f$ on $\cC_\theta$, the capillary $L_p$ dual Minkowski problem considered here is to find a smooth positive function $h: \mathcal{C}_\theta \to \mathbb{R}$ with $\nabla^2h+hI>0$ on $\cC_\theta$ satisfying
\begin{equation}\label{eq:target}
\left\{\begin{aligned}
\det(\nabla^2h+hI)&=f h^{p-1}(h^2+|\nabla h|^2)^{\frac{n+1-q}{2}}, &&\text{in }\operatorname{int}\cC_\theta,\\
\nabla_\mu h&=\cot\theta\,h, &&\text{on }\partial\cC_\theta.
\end{aligned}\right.
\end{equation}
Here $\mu$ is the outward unit conormal of $\partial\cC_\theta$ in $\cC_\theta$, and derivatives use its round metric $g$. The notation $\nabla^2h+hI$ means $\nabla^2h+hg$, represented by $(h_{ij}+h\delta_{ij})$ in an orthonormal frame.

Recent work has established several existence and uniqueness results for \eqref{eq:target}. From the perspective of PDEs, Gao \cite{G26a} proved existence and uniqueness for $p>q$ and $q\leq1$ under an acute contact-angle assumption. The restriction on $q$ arose in the interior $C^2$ estimate. In subsequent work, Gao \cite{G26b} removed this restriction in dimensions $n\geq3$ by constructing a different auxiliary function. Hu, Hu and Yang \cite{HHY26} derived a non-collapsing estimate for strictly convex capillary hypersurfaces using a capillary John ellipsoid theorem. Combining this estimate with gradient estimates, they obtained the required $C^0$ bounds, which then served as a basis for establishing higher-order estimates. For $n\geq2$ and acute contact angles, they proved existence for even data when $1<p<q<\infty$, existence and uniqueness without symmetry assumptions when $p>q$. For even data and $p=q>1$, they solved the normalized eigenvalue problem, with a unique positive eigenvalue and uniqueness of the support function up to positive dilation \cite[Theorem~1.2]{HHY26}. A related extension to the capillary $L_p$ dual Christoffel-Minkowski problem has recently been investigated by Li and Liu \cite{LL26a}.

Geometric flows provide another approach to capillary Minkowski-type problems. Mei, Wang and Weng \cite{MWW25c} studied the capillary Gauss curvature flow. Hu, Hu and Ivaki \cite{HHI25} established long-time existence and convergence results for normalized capillary $L_p$ Minkowski flows with even data for $p>-n-1$, and also treated a flow without symmetry assumptions for $p>n+1$. The purpose of this paper is to extend the geometric flow approach to the capillary $L_p$ dual Minkowski problem for capillary hypersurfaces in the Euclidean half-space.

For $p,q\in\bbR$ and a smooth function $h>0$ with $\nabla^2h+hI>0$ and $\nabla_\mu h=\cot\theta\,h$, set $\rho^2=h^2+|\nabla h|^2$ and
\begin{equation*}
\begin{aligned}
H_p[h]=\int_{\cC_\theta}fh^p\,d\xi,
\qquad Q_q[h]=\int_{\cC_\theta}\rho^{q-n-1}h\det(\nabla^2h+hI)\,d\xi.
\end{aligned}
\end{equation*}

We consider even, smooth, strictly convex capillary hypersurfaces $\Sigma_\tau=X(\mathcal M,\tau)$ with capillary support function $h(\cdot,\tau)$. The initial condition is geometric: $h_0$ is the capillary support function of $\Sigma_0$, and no  parametrization at time zero is prescribed. The normalized flow is
\begin{equation}\label{eq:geometric-flow}
\left\{\begin{aligned}
\partial_\tau X&=-\eta(\tau)\frac{(f\circ\widetilde\nu)\rho^{n+1-q}(X\cdot\nu)^p}{1+\cos\theta(\nu\cdot e)}\cK\widetilde\nu+X, &&\text{in }\mathcal M\times(0,T),\\
\widetilde\nu\cdot e&=0, &&\text{on }\partial\mathcal M\times(0,T),\\
X(\cdot,0)&=X_0 &&\text{in }\mathcal M,
\end{aligned}\right.
\end{equation}
where $\mathcal M$ is a compact manifold with boundary, $\rho=|X|$, $\eta(\tau):=\eta[h(\cdot,\tau)]=\frac{Q_q[h(\cdot,\tau)]}{H_p[h(\cdot,\tau)]}$, $\cK$ is the Gauss curvature of $\Sigma_\tau$,  $e$ denotes the unit vector pointing in the negative direction of the $(n+1)$-th coordinate axis, and $\widetilde\nu$ is the capillary Gauss map defined in \eqref{defn-cG}.

We say that a function $f: \cC_\theta \to \mathbb{R}$ is \emph{even} if
\begin{equation*}
f(-\xi',\xi_{n+1})=f(\xi',\xi_{n+1}) \quad \forall \,\xi=(\xi',\xi_{n+1})\in \cC_\theta.
\end{equation*}
A capillary hypersurface is called even if it is invariant under the horizontal reflection $(x',x_{n+1})\mapsto(-x',x_{n+1})$. Equivalently, its capillary support function defined by \eqref{capillary-support-function} is even. For a solution of \eqref{eq:geometric-flow} with positive curvature speed, differentiating $X\cdot e=0$ at the boundary yields $\widetilde\nu\cdot e=0$, which is the prescribed contact-angle condition.

The long-time existence and convergence results for the flow \eqref{eq:geometric-flow} are summarized as follows.
\begin{thm}\label{thm:main}
Let $n\geq2$,  $p,q\in\bbR$ with $pq\geq0$, and $\theta\in(0,\frac{\pi}{2})$. Let $f$ be an even, smooth, positive function on $\cC_\theta$, and let $\Sigma_0$ be an even, smooth, compact, properly embedded, strictly convex capillary hypersurface in $\overline{\mathbb{R}^{n+1}_+}$. Then there exists an even, smooth, and strictly convex solution $\Sigma_{\tau}$ to the flow \eqref{eq:geometric-flow} for all $\tau > 0$. Moreover, there exists a sequence $\tau_{j}\to\infty$ such that $\Sigma_{\tau_{j}}$ converges smoothly to an even, smooth, strictly convex capillary hypersurface $\Sigma$, whose capillary support function $h$ satisfies
\begin{equation}\label{eq:normalized-target}
\left\{\begin{aligned}
\det(\nabla^2h+h I)&=\eta f h^{p-1}(h^2+|\nabla h |^2)^{\frac{n+1-q}{2}}, &&\text{in }\operatorname{int}\cC_\theta,\\
\nabla_\mu h&=\cot\theta\,h, &&\text{on }\partial\cC_\theta,
\end{aligned}\right.
\end{equation}
where $\eta=Q_q[h]/H_p[h]>0$. In particular, if $p\neq q$, then $\bar h=\eta^{1/(p-q)}h$ solves \eqref{eq:target}.
\end{thm}

The proof of \autoref{thm:main} relies crucially on establishing a priori estimates for the flow \eqref{eq:geometric-flow}. Compared with the $L_p$ dual Minkowski flows, the capillary setting presents additional challenges. In particular, the derivation of the $C^0$ estimate relies on conserved quantities and monotone entropy functionals, as well as a compactness argument based on reflection that is compatible with both the Robin boundary condition and the horizontal symmetry. Moreover, the derivation of the $C^2$ estimates requires carefully constructed auxiliary functions that accommodate the boundary condition while allowing the application of the interior maximum principle. The test function $\log\sigma_1+Ah+B|\nabla(h/\ell)|^2$ and its spatial cancellation are adapted from \cite[Lemma~5.4]{HHY26}. Here the additional time-derivative terms and the nonlocal normalization are controlled along the flow \eqref{eq:geometric-flow}. The main conclusion is long-time existence and subsequential convergence, and stationary existence overlaps with \cite{HHY26} in the parameter ranges covered there.

The paper is organized as follows. In \autoref{Sec2}, we collect some basic facts on the capillary geometry. In \autoref{Sec3}, we derive monotone quantities associated with the flow \eqref{eq:geometric-flow}. The uniform a priori estimates for the flow \eqref{eq:geometric-flow} are established in \autoref{Sec4}. Consequently, the long-time existence and convergence of the flow \eqref{eq:geometric-flow} are established, and the proof of \autoref{thm:main} is completed in \autoref{Sec5}.

\section{Preliminaries}\label{Sec2}

Throughout the paper, $n\geq2$ denotes the dimension of the hypersurfaces. The ambient Euclidean space has dimension $n+1$.

Let $\{E_i\}_{i=1}^{n+1}$ be the standard orthonormal basis of $\bbR^{n+1}$, set $e=-E_{n+1}$, and write $\bbR^{n+1}_+=\{x\in\bbR^{n+1}:x\cdot E_{n+1}>0\}$.  For $r>0$, we define the capillary spherical cap of radius $r$ by
\begin{equation*}
\cC_{\theta,r} =\{\zeta\in\overline{\bbR^{n+1}_+}: |\zeta-r\cos\theta\,e|=r\}.
\end{equation*}
We write $\cC_\theta=\cC_{\theta,1}$ and set
\begin{equation*}
\begin{aligned} \bbS_\theta =\{\omega\in\bbS^n:\omega_{n+1}\geq\cos\theta\}, \qquad \cC_\theta =T(\bbS_\theta),\qquad T(\omega)=\omega+\cos\theta\,e. \end{aligned}
\end{equation*}

Throughout, a strictly convex capillary hypersurface $\Sigma\subset\overline{\bbR^{n+1}_+}$ is smooth, compact and properly embedded, with $\operatorname{int}\Sigma\subset\bbR^{n+1}_+$ and $\partial\Sigma\subset\partial\bbR^{n+1}_+$. It meets the supporting hyperplane at angle $\theta\in(0,\pi/2)$ and, together with its flat side, bounds a capillary convex body $\widehat\Sigma$. Let $\nu$ be its outer unit normal. The capillary Gauss map $\widetilde\nu:\Sigma\to\cC_\theta$ is defined by
\begin{equation}\label{defn-cG}
\widetilde\nu=\nu+\cos\theta\,e. 
\end{equation}
Let $X=\widetilde\nu^{-1}$. The capillary support function $h_\Sigma:\cC_\theta\to\mathbb R$ is defined by
\begin{equation}\label{capillary-support-function}
\begin{aligned} h(\xi):=h_\Sigma(\xi)=X(\xi)\cdot\nu(X(\xi)). \end{aligned}
\end{equation}
It satisfies the Robin boundary condition
\begin{equation}\label{eq:Robin}
\nabla_\mu h=\cot\theta\,h.
\end{equation}
The capillary support function of  $\cC_\theta$ is
\begin{equation*}
\begin{aligned} \ell(\xi)=1+\cos\theta\,(\nu\cdot e) =\sin^2\theta+\cos\theta\,(\xi\cdot e). \end{aligned}
\end{equation*}
In addition,
\begin{equation}\label{eq:ell}
1-\cos\theta\leq\ell\leq\sin^2\theta, \qquad \nabla^2\ell+\ell I=I, \qquad \nabla_\mu \ell=\cot\theta\,\ell.
\end{equation}

Let $\hat h_{\widehat\Sigma}:\bbS^n\to\bbR$ be the standard support function of the capillary convex body $\widehat\Sigma$, defined by
\[
\hat{h}_{\widehat\Sigma}(\omega) := \max_{x \in \widehat{\Sigma}} \omega \cdot x, \quad \forall\, \omega \in \mathbb{S}^n.
\]
By these definitions, we have
\begin{equation}\label{eq:ordinary-capillary-support}
\nu(X(\xi))=T^{-1}(\xi)=\xi-\cos\theta\,e, \qquad h_\Sigma(\xi)=\hat h_{\widehat\Sigma}\bigl(T^{-1}(\xi)\bigr).
\end{equation}

In a local orthonormal frame $\{e_i\}^n_{i=1}$
on $\cC_{\theta}$ (with respect to the standard
metric of $\cC_{\theta}$),  write
\begin{equation*}
b_{ij}=h_{ij}+h\delta_{ij},\qquad \cK=(\det(b_{ij}))^{-1},\qquad \sigma_1=\tr(b_{ij}).
\end{equation*}
Throughout, $\tau$ denotes normalized time and $t$ denotes time for the unnormalized short-time flow. We reserve $u=h/\ell$, and use $\omega=T^{-1}(\xi)$ for the ordinary normal direction and $v=X/|X|$ for the radial direction. Spatial subscripts denote covariant derivatives in a local orthonormal frame, with spatial indices ranging from $1$ to $n$. At the boundary, tangential indices $\alpha,\beta$ range from $1$ to $n-1$, and $e_n=\mu$. Repeated spatial indices are summed, and $(b^{ij})=(b_{ij})^{-1}$. Strict convexity of a support function means $b[h]>0$. The eigenvalues of $b[h]$ are the principal radii of the associated hypersurface. By \cite[Proposition~2.6]{MWWX25}, every positive $h$ satisfying $b[h]>0$ and the Robin boundary condition \eqref{eq:Robin} uniquely determines a smooth strictly convex capillary body.

For an even capillary convex body, its flat side $D_\Sigma:=\widehat\Sigma\cap\{x_{n+1}=0\} \subset\partial\bbR^{n+1}_+$ is a full-dimensional origin-symmetric convex body in the support hyperplane. Thus $0\in\operatorname{relint}D_\Sigma$. Every outer unit normal to $\Sigma$ has vertical component at least $\cos\theta>0$.  $\Sigma$ can be represented as the graph of a function $g:D_\Sigma\to[0,\infty)$ that is positive in $\operatorname{relint}D_\Sigma$. By convexity, $g$ is concave.
Set $a=g(0)>0$. Since $(0,a)\in\Sigma$, for every $\omega\in\bbS_\theta$, we obtain
\begin{equation*}
h(T(\omega))=\hat h_{\widehat\Sigma}(\omega)\geq (0,a)\cdot \omega=a \omega_{n+1}\geq a\cos\theta>0.
\end{equation*}
In particular, the initial capillary support function in \autoref{thm:main} is positive. All capillary support functions used in the flow \eqref{eq:geometric-flow} and variational arguments below are positive and satisfy $b[h]>0$.

The inverse Gauss parametrization satisfies
\begin{equation*}
\begin{aligned}
X(\xi)=\nabla h(\xi)+h(\xi)\nu(X(\xi)), \qquad 
\rho(\xi):=|X(\xi)|=\bigl(h^2+|\nabla h|^2\bigr)^{1/2}.
\end{aligned}
\end{equation*}
Direct differentiation on $\cC_\theta$ yields
\begin{equation}\label{eq:position-length-derivatives}
\rho_i=\frac{b_{ij}h_j}{\rho}, \qquad \left.\frac{\rho_\tau}{\rho}\right|_{\xi} =\frac{hh_\tau+\langle\nabla h,\nabla h_\tau\rangle} {h^2+|\nabla h|^2}.
\end{equation}

For a positive, strictly convex capillary support function $h$, the origin lies in the relative interior of the flat side: its support values in horizontal directions are $h/\sin\theta>0$ at the corresponding boundary normals. Write $\bbS^n_+=\{v\in\bbS^n:v_{n+1}>0\}$ and $\overline{\bbS^n_+}=\{v\in\bbS^n:v_{n+1}\geq0\}$. Radial projection maps $\Sigma$ diffeomorphically onto $\overline{\bbS^n_+}$ as a manifold with boundary. Let $r\in C^\infty(\overline{\bbS^n_+})$ be its positive radial function, so $X=r(v)v$ and $r(v)=\rho(\xi)$ at corresponding points. Its Jacobian is
\begin{equation*}
dv=\frac{X\cdot\nu}{|X|^{n+1}}\,dA
   =h\rho^{-n-1}\det b\,d\xi.
\end{equation*}
Consequently, for every $q\in\bbR$,
\begin{equation*}
Q_q[h]=\int_{\bbS^n_+}r(v)^q\,dv.
\end{equation*}
The radial parametrization is smooth up to the equator, since $v\cdot\nu=h/r>0$. Integrals may be taken over $\bbS^n_+$ because its boundary has zero spherical measure.

For later use, set $u=h/\ell$. Then we have
\begin{equation}\label{eq:u-boundary}
u_\mu=0, \qquad u_{\alpha\mu}=-\cot\theta\,u_\alpha, \qquad 1\leq\alpha\leq n-1,
\end{equation}
and
\begin{equation}\label{eq:b-u}
b_{ij}=\ell u_{ij}+\ell_i u_j+\ell_j u_i+u\delta_{ij}.
\end{equation}
In particular, \eqref{eq:ell}, \eqref{eq:u-boundary}, and \eqref{eq:b-u} imply
\begin{equation*}
b_{\alpha\mu}=0,\qquad 1\leq\alpha\leq n-1.
\end{equation*}

\section{Monotone Quantities}\label{Sec3}

By the definition of the capillary support function and the normalized flow \eqref{eq:geometric-flow}, we have
\begin{equation}\label{eq:scalar-flow}
\left\{\begin{aligned}
\partial_\tau h&=-\eta f\rho^{n+1-q}h^p\cK+h, &&\text{in }\operatorname{int}\cC_\theta\times(0,T),\\
\nabla_\mu h&=\cot\theta\,h, &&\text{on }\partial\cC_\theta\times(0,T),\\
h(\cdot,0)&=h_0, &&\text{in }\cC_\theta.
\end{aligned}\right.
\end{equation}
When $q=n+1$, the identity $Q_{n+1}=(n+1)V(\widehat{\Sigma_\tau})$ shows that \eqref{eq:scalar-flow} reduces to the volume-preserving capillary $L_p$ Minkowski flow (cf. \cite{HHI25}).

\begin{prop}\label{prop:support-variation}
Suppose $q\in \mathbb{R}$. Let $\{h(\cdot,s)\}$ be a smooth family of positive strictly convex functions satisfying $\nabla_{\mu}h(\cdot,s)=\cot\theta\, h(\cdot,s)$. Put $\varphi=\partial_s h(\cdot,s)$ and $\rho(\cdot,s)=(h(\cdot,s)^2+|\nabla h(\cdot,s)|^2)^{1/2}$. Then
\begin{equation}\label{eq:support-first-variation}
\frac d{ds}Q_q[h(\cdot,s)] =q\int_{\cC_\theta} \varphi\,\rho(\cdot,s)^{q-n-1}\det b[h(\cdot,s)]\,d\xi.
\end{equation}
In particular, $Q_0[h(\cdot,s)]$ is independent of $s$.
\end{prop}

\begin{proof}
Use $X(\xi,s)=r(v,s)v$ with $\xi=\xi(v,s)$. At fixed Gauss coordinate, $\partial_sX\cdot\nu=\varphi$. Varying $\xi$ only adds a tangent vector. Thus, at a fixed radial direction $v$,
\begin{equation*}
r_s(v,s)\,v\cdot\nu=\varphi(\xi(v,s),s),
\qquad
v\cdot\nu=\frac{h(\xi(v,s),s)}{r(v,s)}.
\end{equation*}
Therefore $r_s/r=\varphi/h$, with the right-hand side evaluated at the corresponding Gauss coordinate. Using the radial formulas in \autoref{Sec2}, we obtain
\begin{equation*}
\begin{split}
\frac d{ds}Q_q[h(\cdot,s)]
=q\int_{\bbS^n_+}r^q\frac{\varphi}{h}\,dv=q\int_{\cC_\theta}\varphi\rho^{q-n-1}\det b\,d\xi.
\end{split}
\end{equation*}
Differentiation under the integral is justified by smoothness and the positive lower bound for $r$ on the closed hemisphere, locally uniformly in $s$. This proves \eqref{eq:support-first-variation}, including $q=0$.
\end{proof}

The radial formula gives $Q_0[h]=|\bbS^n_+|=|\bbS^n|/2$. Define
\begin{align*}
\mathcal N_q[h]&=\left(\frac{Q_q[h]}{Q_0}\right)^{1/q},\qquad q\neq0,\\
\log\mathcal N_0[h]&=\frac1{Q_0}\left.\partial_qQ_q[h]\right|_{q=0}
=\frac1{Q_0}\int_{\bbS^n_+}\log r(v)\,dv.
\end{align*}

\begin{prop}\label{prop:preserve}
For $q\neq0$, $Q_q[h(\cdot,\tau)]$ is constant along the flow \eqref{eq:scalar-flow}.  For $q=0$, $\mathcal N_0[h(\cdot,\tau)]$ is constant.
\end{prop}

\begin{proof}
For $q\neq0$, \eqref{eq:scalar-flow} and \autoref{prop:support-variation} imply
\begin{equation*}
\frac d{d\tau}Q_q[h(\cdot,\tau)] =q\int_{\cC_\theta} h_\tau \rho^{q-n-1}\det b\,d\xi =q(-\eta H_p+Q_q)=0.
\end{equation*}
For $q=0$, differentiating \eqref{eq:support-first-variation} with respect to the exponent $q$ and then evaluating the resulting identity at $q=0$, yields
\begin{equation*}
\frac d{d\tau}\log\mathcal N_0[h(\cdot,\tau)] =\frac1{Q_0}\int_{\cC_\theta} h_\tau \rho^{-n-1}\det b\,d\xi =\frac{-\eta H_p+Q_0}{Q_0}=0.
\end{equation*}
\end{proof}

For a positive function $g$ on $\cC_\theta$, define the weighted power means
\begin{align*}
\|g\|_{f,s}&=\left(\frac{\int_{\cC_\theta}fg^s\,d\xi}{\int_{\cC_\theta}f\,d\xi}\right)^{1/s}
\quad(s\neq0),\\
\|g\|_{f,0}&=\exp\left(\frac{\int_{\cC_\theta}f\log g\,d\xi}{\int_{\cC_\theta}f\,d\xi}\right).
\end{align*}
Define the functional
\begin{equation*}
\Phi_{p,q}(\Sigma_\tau)=\log\|h(\cdot,\tau)\|_{f,p}
-\log\mathcal N_q[h(\cdot,\tau)].
\end{equation*}

\begin{prop}\label{prop:entropy}
For all $p,q\in\bbR$, along the flow \eqref{eq:scalar-flow}, we obtain
\begin{equation}\label{eq:dissipation}
\frac d{d\tau}\Phi_{p,q}(\Sigma_\tau)=-\int_{\cC_\theta} \frac{(h_\tau)^2}{Q_qh \rho^{n+1-q}\cK}\,d\xi\leq0.
\end{equation}
Equality holds if and only if $\Sigma_\tau$ satisfies \eqref{eq:normalized-target}.
\end{prop}

\begin{proof}
For $p\neq0$, we have
\[
\partial_\tau\log\|h\|_{f,p}=H_p^{-1}\int_{\cC_{\theta}}fh^{p-1}h_\tau\,d\xi.
\]
For $p=0$, the same expression holds with $H_0=\int_{\cC_\theta}f \,d\xi$. Analogously, we have
\[
\partial_\tau\log\mathcal N_q=Q_q^{-1}\int_{\cC_{\theta}} h_\tau\rho^{q-n-1}\det b\,d\xi
\]
for every $q$. Thus
\begin{equation}\label{phi}
\begin{aligned} 
\frac{d}{d\tau}\Phi_{p,q}(\Sigma_\tau) &=\int_{\cC_\theta}h_\tau \left(\frac{fh^{p-1}}{H_p} -\frac{\rho^{q-n-1}}{Q_q\cK}\right)d\xi. 
\end{aligned}
\end{equation}
Since $h-h_\tau=\eta f\rho^{n+1-q}h^p\cK$ and $\eta=Q_q/H_p$,
\begin{equation}\label{fhp}
\frac{fh^{p-1}}{H_p}-\frac{\rho^{q-n-1}}{Q_q\cK}=-\frac{h_\tau}{Q_qh\rho^{n+1-q}\cK}.
\end{equation}
Substituting \eqref{fhp} into \eqref{phi}, we prove \eqref{eq:dissipation}. Equality is equivalent to $h_\tau\equiv0$, and hence to the normalized stationary equation \eqref{eq:normalized-target}.
\end{proof}

\section{A Priori Estimates}\label{Sec4}
 Given a capillary convex hypersurface $\Sigma$, let $\widehat\Sigma$ be the capillary convex body enclosed by $\Sigma$ and  the support hyperplane $\partial \mathbb{R}^{n+1}_{+}$. Set $\mathcal R(x',x_{n+1})=(x',-x_{n+1})$ and $K=\widehat\Sigma\cup\mathcal R\widehat\Sigma$. Its Minkowski functional is
\begin{equation*}
g_K(x)=\inf\{\lambda>0:x\in\lambda K\}.
\end{equation*}
For the evolving capillary convex body $\widehat{\Sigma_\tau}$, denote the corresponding doubled body by $K_\tau$. For an origin-symmetric convex body containing the origin, write
\begin{equation*}
R(K)=\max_{x\in K}|x|, \qquad r(K)=\max\{s>0:sB^{n+1}\subset K\} =\frac1{\max_{\bbS^n}g_K}.
\end{equation*}
Here $B^{n+1}$ is the Euclidean unit ball in $\mathbb R^{n+1}$. For $\xi\in\cC_\theta$, write
\begin{equation*}
\omega=T^{-1}(\xi)=\xi-\cos\theta\,e\in\bbS_\theta.
\end{equation*}
Translation identifies the spherical measures on $\cC_\theta$ and $\bbS_\theta$.

\begin{lemma}\label{lem:doubling}
Let $\widehat\Sigma$ be the body enclosed by a smooth strictly convex capillary hypersurface $\Sigma$ and its flat side, with positive horizontally even support function and the fixed acute contact angle $\theta \in (0,\frac{\pi}{2})$. Then $K$ is an origin-symmetric convex body containing the origin in its interior. Moreover,
\begin{equation}\label{eq:doubling-data}
\hat h_K(T^{-1}(\xi))=h(\xi),\quad \xi\in\cC_\theta, \qquad 2Q_q[h]=\int_{\bbS^n}g_K(v)^{-q}\,dv.
\end{equation}
\end{lemma}

\begin{proof}
Let $D_{\Sigma}=\widehat\Sigma\cap\{x_{n+1}=0\}$ be the flat side of $\widehat\Sigma$. The Gauss image of $\Sigma$ is $\bbS_\theta$, so the vertical component of every normal is at least $\cos\theta>0$. Thus the vertical projection is locally one-to-one on $\Sigma$. Strict convexity makes it globally one-to-one in the interior, and the contact boundary projects onto $\partial D_{\Sigma}$. Hence $\Sigma$ is the graph of a concave function $g:D_{\Sigma}\to[0,\infty)$, and convexity shows that every vertical section of $\widehat\Sigma$ is the interval joining the flat side to this graph. Therefore
\begin{equation*}
\widehat\Sigma=\{(x',s):x'\in D_{\Sigma},\ 0\leq s\leq g(x')\}.
\end{equation*}
Consequently, $K=\{(x',s):x'\in D_{\Sigma},\ |s|\leq g(x')\}$. If $|s_i|\leq g(x_i')$, concavity of $g$ implies
\begin{equation*}
|\lambda s_1+(1-\lambda)s_2| \leq\lambda g(x_1')+(1-\lambda)g(x_2') \leq g(\lambda x_1'+(1-\lambda)x_2'),
\end{equation*}
which proves convexity. Horizontal reflection preserves $\bbS_\theta$, and \eqref{eq:ordinary-capillary-support} identifies the horizontal evenness of $h$ with the invariance of $\widehat\Sigma$ under $(x',s)\mapsto(-x',s)$. Hence $D_{\Sigma}=-D_{\Sigma}$ and $g(-x')=g(x')$. Reflection across $\{x_{n+1}=0\}$ then implies $K=-K$.

At a point $X(\xi)\in\Sigma$, $T^{-1}(\xi)$ is its outer unit normal. Since its vertical component is positive,  one sees  $\omega'\cdot y\leq h(\xi)$ for $y\in D_{\Sigma}$. Therefore, for $(y,-s)$ in the reflected half $\mathcal{R}\widehat\Sigma$,
\begin{equation*}
\omega\cdot(y,-s)=\omega'\cdot y-\omega_{n+1}s\leq h(\xi).
\end{equation*}
The same supporting hyperplane therefore supports $\mathcal{R}\widehat\Sigma$, and
\begin{equation*}
\hat h_K(T^{-1}(\xi)) =\langle X,T^{-1}(\xi)\rangle=h(\xi).
\end{equation*}

For the second identity, the radial parametrization in \autoref{Sec2} gives $r(v)=g_K(v)^{-1}$ for $v\in\bbS^n_+$. Hence
\begin{equation}\label{eq:radial-change-of-variables}
\int_{\bbS^n_+}g_K(v)^{-q}\,dv
=\int_{\cC_\theta}h\rho^{q-n-1}\det b\,d\xi
=Q_q[h],
\end{equation}
where $g_K(v)^{-1}=|X|$ at the corresponding point. Reflection produces the same integral on the lower hemisphere. At a ridge point $y=(y',0)\in\partial D_{\Sigma}\times\{0\}$, the radial direction is $(y'/|y'|,0)$. Hence the radial image of the entire ridge lies in the equator and has zero spherical measure. Note that the supporting normals of the ridge can have positive spherical measure. The integral in \eqref{eq:radial-change-of-variables} is over radial directions, so this normal image does not create an additional term. Adding the two hemisphere integrals yields \eqref{eq:doubling-data}. The doubled body is convex, with no smoothness assumption at its ridge.
\end{proof}

\begin{lemma}\label{lem:symmetric-compactness}
Let $\Omega\subset\bbS^n$ be measurable and contain a nonempty open set, and let $w$ be measurable with $0<c\leq w\leq C$. There exist constants $a_s>0$ for $s>0$ and $c_0\in\mathbb R$, depending only on $n,\Omega,c,C$ and, for $a_s$, on $s$, such that the following conclusions hold.

(i) For any $a\in\bbS^n$,
\begin{equation}\label{eq:angular-lower}
\int_\Omega w(\omega)|a\cdot \omega|^s\,d\omega\geq a_s \ (s>0), \qquad \int_\Omega w(\omega)\log|a\cdot \omega|\,d\omega\geq c_0.
\end{equation}

(ii) If \eqref{eq:angular-lower} holds, then for every origin-symmetric convex body $K$ with outer radius $R(K)$,
\begin{equation}\label{RK}
\begin{gathered}
\int_\Omega w(\omega)\hat h_K(\omega)^s\,d\omega\geq a_sR(K)^s \quad(s>0),\\
\int_\Omega w(\omega)\log\hat h_K(\omega)\,d\omega\geq\left(\int_\Omega w(\omega)\,d\omega\right)\log R(K)+c_0.
\end{gathered}
\end{equation}

(iii) If $g_K$ is the Minkowski functional of an origin symmetric convex body  $K$, then
\begin{equation}\label{eq:gauge-angular-lower}
\begin{gathered}
\int_\Omega w g_K^s\geq a_s\big(\max_{\bbS^n}g_K\big)^s \quad(s>0),\\
\int_\Omega w\log g_K\geq c_0+\left(\int_\Omega w(v)\,dv\right)\log\max_{\bbS^n}g_K.
\end{gathered}
\end{equation}

(iv) Let $\{K_j\}$ be origin-symmetric convex bodies whose inradii tend to zero. Then, after passing to a subsequence,
\begin{equation}\label{eq:collapse-ae}
g_{K_j}(v)^{-1}\longrightarrow0 \quad\text{for almost every }v\in\bbS^n.
\end{equation}
\end{lemma}

\begin{proof}
Let $a_0\in\mathbb S^n$ be such that $R(K)a_0\in\partial K$. By the definition of the standard support function, since $K$ is origin symmetric, we have
\begin{equation*}
\hat h_K(\omega)\geq R(K)|a_0\cdot \omega|.
\end{equation*}
Then \eqref{eq:angular-lower}  implies
\begin{equation*}
\begin{aligned}
\int_{\Omega}{w\hat{h}_K^s}\,d\omega\geq a_sR(K)^s,\qquad
\int_{\Omega}{w\log{\hat{h}_K}}\,d\omega \geq \left(\int_\Omega w(\omega)\,d\omega\right)\log R(K)+c_0.
\end{aligned}
\end{equation*}
This proves \eqref{RK}.

Now we prove \eqref{eq:angular-lower}. For the first inequality, given $s>0$ and $a\in \bbS^n$, the integral $\int_\Omega |a\cdot \omega|^s\,d\omega$ is positive, since $\Omega$ contains a nonempty open set and $a^\perp\cap\mathbb S^n$ has empty interior. The map $a\mapsto\int_\Omega |a\cdot \omega|^s\,d\omega$ is continuous, so compactness of $\mathbb{S}^n$ gives $\int_\Omega |a\cdot \omega|^s\,d\omega\geq c^*>0$, and it follows that $\int_\Omega w(\omega)|a\cdot \omega|^s\,d\omega\geq cc^*$. For the logarithmic inequality, fix $a\in\bbS^n$ and use $\log|a\cdot \omega|\leq0$ and $w\leq C$ to obtain
\begin{equation*}
\int_\Omega w\log|a\cdot \omega|\,d\omega
\geq C\int_{\bbS^n}\log|a\cdot \omega|\,d\omega\geq \widetilde{C},
\end{equation*}
where we used the fact that the  integral $\int_{\bbS^n} \log|a\cdot \omega|\, d\omega$ is finite and independent of $a$ by rotational invariance.

Let $K^\circ=\{y\in\bbR^{n+1}:x\cdot y\leq1\ \text{for all }x\in K\}$ be the polar body. The standard polar identities are
\begin{equation*}
g_K=\hat h_{K^\circ},\qquad
R(K^\circ)=\max_{\bbS^n}g_K=\frac1{r(K)}.
\end{equation*}
Since $K^\circ$ is also origin symmetric, (iii) follows by applying (ii) to $K^\circ$.

For (iv), set $M_j=R(K_j^\circ)=1/r(K_j)\to\infty$ and choose $a_j\in\bbS^n$ with $M_ja_j\in K_j^\circ$. After passing to a subsequence, $a_j\to a\in\bbS^n$. The origin symmetry of $K_j^\circ$ implies
\begin{equation*}
g_{K_j}(v)=\hat h_{K_j^\circ}(v)
\geq M_j|a_j\cdot v|\longrightarrow\infty
\quad\text{if }a\cdot v\neq0.
\end{equation*}
Thus $g_{K_j}(v)^{-1}\to0$ outside the null set $a^\perp\cap\bbS^n$, proving \eqref{eq:collapse-ae}.
\end{proof}

We first deduce the $C^0$ estimate.
\begin{lemma}\label{lem:C0}
Let $pq\geq 0$ and $\theta \in (0,\frac{\pi}{2})$. Suppose $h(\cdot, \tau)$ is a positive, even, smooth strictly convex solution to the flow \eqref{eq:scalar-flow}. Then there exists a positive constant $C$, depending only on $n, p, q$, $\theta$, $\min_{\mathcal{C}_\theta}f, \max_{\mathcal{C}_\theta}f$ and $\Sigma_0$, such that
\begin{equation}\label{eq:C0}
C^{-1}\leq h(\xi,\tau)\leq C, \quad  \forall (\xi,\tau) \in \mathcal{C}_\theta \times [0,T).
\end{equation}
\end{lemma}

\begin{proof}
We identify $f$ with its pullback $w=f\circ T$ to $\bbS_\theta$, which has the same positive upper and lower bounds. At $q=0$, \eqref{eq:doubling-data}  yields $2Q_0[h]=|\bbS^n|$ and
\begin{equation}\label{eq:N0-gauge}
\log\mathcal N_0[h] =-\frac1{|\bbS^n|}\int_{\bbS^n}\log g_{K_\tau}\,dv.
\end{equation}
By \autoref{prop:preserve} and \autoref{prop:entropy}, we obtain
\begin{equation}\label{upp-h-f}
\log\|h(\cdot,\tau)\|_{f,p} \leq \Phi_{p,q}(\Sigma_0)+\log\mathcal N_q[h_0].
\end{equation}
This gives a positive constant $C$, independent of $\tau$, such that
\begin{align}
\int_{\cC_\theta} f h(\xi,\tau)^p\,d\xi
&\leq C, \qquad p>0,\label{a1}\\
\int_{\cC_\theta} f\log h(\xi,\tau)\,d\xi
&\leq C, \qquad p=0,\label{a2}\\
\int_{\cC_\theta} f h(\xi,\tau)^p\,d\xi
&\geq C, \label{a3}\qquad p<0.
\end{align}

We next consider the cases $q>0$ and $q\leq0$. Suppose first that $q>0$. Since $pq\geq0$, we have $p\geq0$. By \eqref{a1} and \eqref{a2}, we obtain
\begin{align*}
\int_{\bbS_\theta} w\,\hat h_{K_\tau}^{\,p}\,d\omega
&\leq C, \qquad p>0,\\
\int_{\bbS_\theta} w\log\hat h_{K_\tau}\,d\omega
&\leq C,
\qquad p=0.
\end{align*}
By \eqref{RK}, there exists a constant $C_0>0$, independent of $\tau$, such that
\begin{equation}\label{R0}
R(K_\tau)\leq C_0.
\end{equation}

We now prove that the inradius has a uniform positive lower bound. Suppose otherwise. Then there exists a sequence of times $\tau_j$ such that
\[
r(K_{\tau_j})\longrightarrow0.
\]
By \eqref{eq:collapse-ae}, after passing to a subsequence,
\[
g_{K_{\tau_j}}(v)^{-1}\longrightarrow0
\qquad\text{for almost every }v\in\bbS^n.
\]
Since $q>0$, it follows that $g_{K_{\tau_j}}(v)^{-q}\to0$ for almost every $v$. Moreover, from \eqref{R0}, we obtain
\[
0\leq g_{K_{\tau_j}}(v)^{-q}
\leq C_0^q.
\]
The constant $C_0^q$ is integrable on $\bbS^n$. Therefore, the dominated convergence theorem yields
\[
\int_{\bbS^n}g_{K_{\tau_j}}(v)^{-q}\,dv
\longrightarrow0.
\]

On the other hand, \eqref{eq:doubling-data} and \autoref{prop:preserve} imply that, for every $j$,
\[
\int_{\bbS^n}g_{K_{\tau_j}}(v)^{-q}\,dv
=2Q_q[h(\cdot,\tau_j)]
=2Q_q[h_0]>0.
\]
This contradiction shows that there exists a constant $r_0>0$, independent of $\tau$, such that
\[
r(K_\tau)\geq r_0.
\]

Suppose next that $q\leq0$. Set
\[
g_\tau=g_{K_\tau},
\qquad
M_\tau=\max_{\bbS^n}g_\tau
=\frac1{r(K_\tau)}.
\]
If $q<0$, put $s=-q>0$. Applying \eqref{eq:gauge-angular-lower} with $\Omega=\bbS^n$ and $w\equiv1$, and using \autoref{prop:preserve}, we obtain
\[
a_sM_\tau^s
\leq\int_{\bbS^n}g_\tau^s\,dv
=2Q_q[h_0].
\]
Consequently,
\[
M_\tau\leq
\left(\frac{2Q_q[h_0]}{a_s}\right)^{1/s},
\qquad
r(K_\tau)\geq
\left(\frac{a_s}{2Q_q[h_0]}\right)^{1/s}>0.
\]
If $q=0$, \autoref{prop:preserve} and \eqref{eq:N0-gauge} give
\begin{equation}\label{G1}
\int_{\bbS^n}\log g_\tau\,dv
=-|\bbS^n|\log\mathcal N_0[h_0].
\end{equation}
On the other hand,  using \eqref{eq:gauge-angular-lower}, applied with $\Omega=\bbS^n$ and $w\equiv1$, yields
\begin{equation}\label{G2}
c_0+|\bbS^n|\log M_\tau
\leq\int_{\bbS^n}\log g_\tau\,dv.
\end{equation}
Combining \eqref{G1} and \eqref{G2}, we find
\[
\log M_\tau
\leq-\log\mathcal N_0[h_0]
-\frac{c_0}{|\bbS^n|}.
\]
Thus
\[
r(K_\tau)
=\frac1{M_\tau}
\geq
\mathcal N_0[h_0]
\exp\left(\frac{c_0}{|\bbS^n|}\right)>0.
\]
In both cases, there exists a constant $r_0>0$, independent of $\tau$, such that
\begin{equation}\label{ro}
r(K_\tau)\geq r_0.
\end{equation}

When $p\geq0$, irrespective of the sign of $q$, \eqref{a1}--\eqref{a2} and \eqref{RK} also give $R(K_\tau)\leq C$.

The remaining step is to prove a uniform upper bound for the outer radius when $p<0$. By \eqref{ro}, there exists a constant $c>0$, independent of $\tau$, such that
\[
cB^{n+1}\subset K_\tau.
\]
From \eqref{a3}, we conclude that there exists  some $C^*>0$, independent of $\tau$, such that
\begin{equation}\label{WV}
\int_{\bbS_\theta} w(\omega)\hat h_{K_\tau}(\omega)^p\,d\omega
\geq C^{*}.
\end{equation}

Suppose, for contradiction, that the outer radii are not uniformly bounded. Then there exists a sequence of times $\tau_j$ such that
\[
R_j:=R(K_{\tau_j})\longrightarrow\infty.
\]
Choose $x_j\in\partial K_{\tau_j}$ with $|x_j|=R_j$, and set $ a_j=\frac{x_j}{R_j}\in\bbS^n$. By compactness of $\bbS^n$, after passing to a subsequence, we may assume that $ a_j\longrightarrow a \in\bbS^n.$

Since $cB^{n+1}\subset K_{\tau_j}$, we have
\begin{equation}\label{hc}
\hat h_{K_{\tau_j}}(\omega)
\geq c
\qquad\text{for every }\omega\in\bbS^n.
\end{equation}
Moreover, since $K_{\tau_j}$ is origin-symmetric,  we obtain
\[
\hat h_{K_{\tau_j}}(\omega)
\geq R_j|a_j\cdot \omega|.
\]
In view of the above, one sees
\[
\hat h_{K_{\tau_j}}(\omega)
\geq\max\{c,R_j|a_j\cdot \omega|\}.
\]

Define
\[
E=\{\omega\in\bbS_\theta:a \cdot \omega=0\}.
\]
Then $E$ is contained in the intersection of $\bbS^n$ with a hyperplane through the origin. Thus $E$ has zero spherical measure.

Now fix $\omega\in\bbS_\theta\setminus E$, and set
\[
\delta_{\omega}=\frac12|a \cdot \omega|>0.
\]
Since $a_j\to a $, there exists an index $j_0(\omega)$ such that
\[
|a_j\cdot \omega|\geq\delta_{\omega}
\qquad\text{for all }j\geq j_0(\omega).
\]
Consequently,
\[
\hat h_{K_{\tau_j}}(\omega)
\geq R_j\delta_{\omega}\longrightarrow\infty.
\]

Since $p<0$,
\[
\hat h_{K_{\tau_j}}(\omega)^p\longrightarrow0, \qquad\text{as}\,\,j\to\infty,\,
\text{for almost every }\omega\in\bbS_\theta.
\]
On the other hand, for $p<0$, \eqref{hc} gives
\[
0\leq w(\omega)\hat h_{K_{\tau_j}}(\omega)^p
\leq\|w\|_{L^\infty(\bbS_\theta)}c^p.
\]
Dominated convergence gives
\[
\int_{\bbS_\theta}
w(\omega)\hat h_{K_{\tau_j}}(\omega)^p\,d\omega
\longrightarrow0.
\]
This contradicts \eqref{WV}. Consequently, there exists a constant $C>0$, independent of $\tau$, such that
\[
R(K_\tau)\leq C.
\]

In all cases, the inner and outer radii of $K_\tau$ satisfy $r(K_\tau)\geq c$ and $R(K_\tau)\leq C$. Since $r(K_\tau)\leq\hat h_{K_\tau}\leq R(K_\tau)$ and $\hat h_{K_\tau}\circ T^{-1}=h$, this proves the desired $C^0$ estimate \eqref{eq:C0}.

\end{proof}

The $C^0$ estimate also gives the following bounds.
\begin{lemma}\label{lem:C1}
Let $pq\geq 0$ and $\theta \in (0,\frac{\pi}{2})$. Suppose $h(\cdot, \tau)$ is a positive, even, smooth strictly convex solution to the flow \eqref{eq:scalar-flow}. Then there exists a positive constant $C$, depending only on $n, p, q$, $\theta$, $\min_{\mathcal{C}_\theta}f, \max_{\mathcal{C}_\theta}f$ and $\Sigma_0$, such that
\begin{equation}\label{eq:C1}
|\nabla h(\cdot,\tau)|\leq C, \qquad C^{-1}\leq \rho(\cdot,\tau)\leq C,
\end{equation}
and
\begin{equation}\label{eta}
C^{-1}\leq\eta(\tau)\leq C.
\end{equation}
The constant has the same dependence as in \autoref{lem:C0}.
\end{lemma}

\begin{proof}
The proof of \autoref{lem:C0} gives $|X|\leq R(K_\tau)\leq C$. Thus $|\nabla h|\leq\rho=|X|\leq C$ and $\rho\geq h\geq c$. The radial formula for $Q_q$ and the bounds for $h$ give $c\leq Q_q,H_p\leq C$, hence $c\leq\eta=Q_q/H_p\leq C$.
\end{proof}

Next, our goal is to establish the $C^2$ estimates. To this end, we first derive upper and lower bounds for the Gauss curvature.

\begin{lemma}\label{prop:Gauss-bounds}
Let $pq\geq 0$ and $\theta \in (0,\frac{\pi}{2})$. Suppose $h(\cdot, \tau)$ is a positive, even, smooth strictly convex solution to the flow \eqref{eq:scalar-flow}. Then there exist positive constants $c$ and $C$, depending only on $n, p, q$, $\theta$, $\min_{\mathcal{C}_\theta}f, \max_{\mathcal{C}_\theta}f$ and $\Sigma_0$, such that
\begin{equation}\label{Bound:K}
c\leq\cK(\xi,\tau)\leq C.
\end{equation}
\end{lemma}

\begin{proof}
\emph{Step 1:} Prove $\cK\leq C$.

Let $c_h>0$ be the time-independent lower bound for $h$ in \autoref{lem:C0}, and fix
\begin{equation*}
\varepsilon_0=c_h/2>0,
\end{equation*}
and set
\begin{equation}\label{eq:Qplus}
Q_{+}=\frac{\eta f \rho^{n+1-q}h^p\cK-h}{h-\varepsilon_0}=\frac{-h_\tau}{h-\varepsilon_0},
\end{equation}
where $h_\tau=\partial_\tau h$. The $C^0$ and $C^1$ estimates give $c\cK-C\leq Q_+\leq C\cK$. It therefore suffices to bound $Q_+$ above.

Fix an arbitrary $T_0\in(0,T)$ and take a maximum of $Q_{+}$ on $\cC_\theta\times[0,T_0]$. If this maximum is non-positive or occurs at $\tau=0$, the bound follows. Otherwise, the chosen maximum point $(\xi_0,\tau_0)$ with $\tau_0>0$ satisfies $\partial_\tau Q_{+}\geq0$, and $Q_{+}(\xi,\tau_0)\leq Q_{+}(\xi_0,\tau_0)$ for every $\xi$. If this positive maximum were attained on $\partial\cC_\theta$, then at that point we would have
\[
\nabla_\mu Q_{+}= -\frac{h_{\tau\mu}}{h-\varepsilon_0} + \frac{h_{\tau} h_{\mu}}{(h - \varepsilon_0)^2} =\frac{\varepsilon_0 \cot\theta\, h_{\tau}}{(h - \varepsilon_0)^2} = -\frac{\varepsilon_0 \cot\theta\, Q_{+}}{h - \varepsilon_0} < 0.
\]
This contradicts $\nabla_\mu Q_{+}\geq0$ at a boundary maximum. Therefore, the chosen maximum $(\xi_0,\tau_0)$ satisfies $\xi_0\in\cC_\theta\setminus\partial\cC_\theta$.

Choose a local orthonormal frame $\{e_i\}_{i=1}^n$ such that $\{b_{ij}\}$ is diagonal at the point $\xi_0$. Then, at the point $\xi_0$, we have
\begin{equation}\begin{split}\label{Up1}
0 = \nabla_i Q_{+} &= \frac{-h_{\tau i}}{h - \varepsilon_0} + \frac{h_{\tau} h_i}{(h - \varepsilon_0)^2},\\
0 \geq \nabla^2_{ii} Q_{+} &= \frac{-h_{\tau ii}}{h - \varepsilon_0} + \frac{2 h_{\tau i} h_i + h_{\tau} h_{ii}}{(h - \varepsilon_0)^2} - \frac{2 h_{\tau} h_i^2}{(h - \varepsilon_0)^3} \\
&= \frac{-h_{\tau ii}}{h - \varepsilon_0} + \frac{h_{\tau} h_{ii}}{(h - \varepsilon_0)^2}.
\end{split}
\end{equation}
This implies that
\begin{equation}\begin{split}
\label{Up3}
-h_{\tau ii} - h_{\tau} &\leq -\frac{h_{\tau} h_{ii}}{h - \varepsilon_0} - h_{\tau} = \frac{-h_{\tau}}{h - \varepsilon_0} [h_{ii} + (h - \varepsilon_0)] = Q_{+} (b_{ii} - \varepsilon_0).
\end{split}
\end{equation}

On the other hand,
\begin{equation*}
\partial_\tau Q_{+}=\frac{h_\tau-h_{\tau\tau}}{h-\varepsilon_0}+Q_{+}+Q_{+}^2.
\end{equation*}
Differentiating \eqref{eq:scalar-flow} in time yields
\begin{equation*}
\frac{h_\tau-h_{\tau\tau}}{h-h_\tau}=\frac{\eta_\tau}{\eta} +(n+1-q)\frac{\rho_\tau}{\rho}+p\frac{h_\tau}{h} +\frac{\cK_\tau}{\cK}.
\end{equation*}
Moreover, by a direct computation,
\begin{equation*}
\frac{h-h_\tau}{h-\varepsilon_0}=Q_{+}+\frac{h}{h-\varepsilon_0},
\end{equation*}
and
\begin{equation}\label{eq:Kt-upper}
\frac{\cK_\tau}{\cK} \leq nQ_{+} \left(1-\varepsilon_0\cK^{\frac{1}{n}}\right).
\end{equation}
Also, we have
\begin{equation}\label{eq:eta-time}
\frac{\eta_\tau}{\eta} =-p\frac{\int_{\cC_\theta}fh^{p-1}h_\tau\,d\xi}{H_p}.
\end{equation}
Since $h-h_\tau>0$ and $h\geq2\varepsilon_0$, we have $Q_{+}\geq-h/(h-\varepsilon_0)\geq-2$. Applying $C^0$, $C^1$ estimates, using $h_\tau=-(h-\varepsilon_0)Q_{+}$ at every spatial point and the global spatial maximum at $\tau_0$, one finds
\begin{equation*}
\left|\frac{\eta_\tau}{\eta}\right|\leq C(1+Q_{+}).
\end{equation*}

At the chosen maximum point, \eqref{eq:Qplus} and \eqref{eq:position-length-derivatives} give
\begin{equation}\label{eq:rho-time-Qplus}
\frac{\rho_\tau}{\rho} =-Q_{+}\left(1-\frac{\varepsilon_0h}{\rho^2}\right).
\end{equation}
Employing \eqref{eq:eta-time}, \eqref{eq:rho-time-Qplus}, and \eqref{eq:Kt-upper},
\begin{equation*}
\begin{aligned}
\partial_\tau Q_{+} \leq{}&\left(Q_{+}+\frac{h}{h-\varepsilon_0}\right)
\left[C(1+Q_{+})+nQ_{+}-n\varepsilon_0Q_{+}\cK^\frac{1}{n}\right]+Q_{+}+Q_{+}^2.
\end{aligned}
\end{equation*}
Since $cQ_{+}\leq\cK\leq C(1+Q_{+})$, it follows that at $(\xi_0,\tau_0)$, we have
\begin{equation*}
\begin{aligned}
0\leq \partial_\tau Q_{+}\leq C(1+Q_{+})^2-cQ_{+}^{2+1/n}.
\end{aligned}
\end{equation*}
The right-hand side is negative when $Q_+$ exceeds a fixed threshold. Hence there is a constant $C_1$, independent of $T_0<T$, such that
\[
Q_{+}(\cdot,\tau) \leq C_1,\qquad 0\leq\tau\leq T_0.
\]
Since $T_0<T$ is arbitrary, the bound holds throughout $[0,T)$. Consequently,
\[
\cK(\cdot,\tau)\leq C_2,
\]
where $C_2$ depends only on the stated data.

\emph{Step 2:} Prove $\cK\geq c$.

We apply the maximum principle to the auxiliary function
\[
Z=f\rho^{n+1-q}\cK, \qquad  Q_{-}=\log Z^{-1}-N\log h,
\]
where $N>\max\{0,p,q-1\}$ is fixed. The $C^0$ and $C^1$ estimates give $|Q_-+\log\cK|\leq C$, so it suffices to bound $Q_-$ above.

\begin{claim}
For this choice of $N$, $\nabla_\mu Q_-<0$ on $\partial\cC_\theta$.
\end{claim}
To prove the claim, fix $\xi_0\in\partial\cC_\theta$. By definition, $h-h_\tau=\eta h^pZ$. Since $h$ and $h_\tau$ satisfy the same Robin condition,
\begin{equation}\label{h-ht-b}
\frac{(h-h_{\tau})_\mu}{h-h_{\tau}}=\cot\theta,
\end{equation}
On the other hand,
\begin{equation}\label{h-ht-b-1}
\frac{(h-h_{\tau})_\mu}{h-h_{\tau}}=p\frac{h_{\mu}}{h}+\frac{Z_{\mu}}{Z}=p\cot\theta+\frac{Z_{\mu}}{Z}.
\end{equation}
Combining \eqref{h-ht-b} with \eqref{h-ht-b-1} gives
\begin{equation}
\frac{Z_{\mu}}{Z}=(1-p)\cot\theta,
\end{equation}
and $N>p-1$ gives
\[
\nabla_\mu  Q_{-} =-\frac{Z_{\mu}}{Z}-N\frac{h_{\mu}}{h}=(p-1-N)\cot\theta<0.
\]
Hence, fix an arbitrary $T_0\in(0,T)$ and take a maximum of $Q_{-}$ on $\cC_\theta\times[0,T_0]$. If it occurs at $\tau=0$, it is bounded by the initial data. Otherwise, any positive-time maximum $(\xi_0,\tau_0)$  of $Q_{-}$ with $\tau_0>0$ satisfies $\xi_0\in\cC_\theta\setminus\partial\cC_\theta$, having $\partial_\tau Q_{-}\geq0$, and $Q_{-}(\xi,\tau_0)\leq Q_{-}(\xi_0,\tau_0)$ for every $\xi$.

At the point $\xi_0$, the condition $\nabla Q_{-}=0$ and $\nabla^2 Q_{-}\leq 0$ give
\begin{equation}\label{eq:Z-gradient}
\frac{Z_i}{Z}=-N\frac{h_i}{h}
\end{equation}
and
\begin{equation}\label{eq-zh}
0\geq (Q_-)_{ii}=-\frac{Z_{ii}}{Z}+\frac{Z_i^2}{Z^2}-N\frac{h_{ii}}{h}+N\frac{h_i^2}{h^2},
\end{equation}
which, together with \eqref{eq:Z-gradient}, implies
\begin{equation}\label{eq:Z-hessian}
\frac{Z_{ii}}Z \geq -N\frac{h_{ii}}h+(N+N^2)\frac{h_i^2}{h^2}.
\end{equation}
If we write
\begin{equation*}
z=\eta h^{p-1}Z=\frac{h-h_\tau}{h},
\end{equation*}
then differentiating $h_\tau=-\eta h^pZ+h$ twice gives
\begin{equation}\label{twice-h}
h_{\tau ii}=h_{ii}-\eta(h^p Z)_{ii}.
\end{equation}
For the term $(h^p Z)_{ii}$, we calculate as follows:
\begin{equation}\label{eq-hpZ}
\begin{aligned}
\frac{(h^p Z)_{ii}}{h^p Z}&=\frac{Z_{ii}}{Z}+p\frac{h_{ii}}{h}+p(p-1)\frac{h_i^2}{h^2}+2p\frac{h_iZ_i}{hZ}\\
&\geq (p-N)\frac{h_{ii}}{h}+(N-p)(N-p+1)\frac{h_i^2}{h^2},
\end{aligned}
\end{equation}
where we used \eqref{eq:Z-gradient} and \eqref{eq:Z-hessian} in the inequality. Combining \eqref{twice-h} with \eqref{eq-hpZ} and noting the fact that $\eta h^p Z=zh$ gives
\begin{equation}\label{eq:htt-bound}
h_{\tau ii}\leq \left[1+(N-p)z\right]h_{ii} -(N-p)(N-p+1)\frac zh h_i^2.
\end{equation}
At the maximum of $Q_-$, \eqref{eq:Z-gradient} and $h_\tau=h(1-z)$ give
\begin{equation}\label{eq:h-tau-i}
h_{\tau i}=h_i-\eta(h^pZ)_i
=\left[1+(N-p)z\right]h_i.
\end{equation}
Moreover,
\begin{equation}\label{eq-rho-tau}
\begin{aligned}
\frac{\rho_\tau}{\rho}&=\frac{hh_{\tau}+\langle\nabla h,\nabla h_\tau\rangle}{\rho^2}\\
 &=\frac{h^2(1-z)+\left[1+(N-p)z\right]|\nabla h|^2}{\rho^2}\\
 &=1+z\frac{(N-p)|\nabla h|^2-h^2}{\rho^2},
\end{aligned}
\end{equation}
where we  used \eqref{eq:h-tau-i} in the second equality and $\rho^2=h^2+|\nabla h|^2$ in the last equality. Therefore,
\begin{equation}\label{eq:Q-cal}
\begin{aligned}
\partial_\tau Q_{-}&=-\frac{Z_{\tau}}{Z}-N\frac{h_{\tau}}{h}\\
&=-\frac{\cK_{\tau}}{\cK}-(n+1-q)\frac{\rho_{\tau}}{\rho}-N\frac{h_{\tau}}{h}\\
&=\sum_i b^{ii}(h_{\tau ii}+h_\tau)-(n+1-q)\frac{\rho_\tau}{\rho}-N\frac{h_\tau}{h}.
\end{aligned}
\end{equation}

By \eqref{eq:htt-bound}, we have
\begin{align*}
h_{\tau ii}+h_{\tau}&\leq\left[1+(N-p)z\right](b_{ii}-h) -(N-p)(N-p+1)\frac zh h_i^2+h(1-z)\\
&=\left[1+(N-p)z\right]b_{ii}-(N-p+1)hz-(N-p)(N-p+1)\frac zh h_i^2
\end{align*}
and
\begin{equation}\label{b-h-t}
\begin{aligned}
\sum_i b^{ii}(h_{\tau ii}+h_\tau)\leq {}& n\left[1+(N-p)z\right]-(N-p+1)hz\sum_i b^{ii}\\
{}&-\frac{(N-p)(N-p+1)z}{h}\sum_{i}b^{ii}h_i^2.
\end{aligned}
\end{equation}
Substituting \eqref{eq-rho-tau} and \eqref{b-h-t} into \eqref{eq:Q-cal} yields
\begin{equation*}
\begin{aligned}
 \partial_\tau Q_{-}\leq{}&q-1-N+z\left[n(N-p)+N\right]+z(n+1-q)\frac{h^2-(N-p)|\nabla h|^2}{\rho^2}\\
&-(N-p+1)z\left(h\sum_i b^{ii}
 +\frac{N-p}{h}\sum_i b^{ii}h_i^2\right).
\end{aligned}
\end{equation*}
The last line is non-positive because $N>p$ and $z>0$. Since $\rho^2=h^2+|\nabla h|^2$,
\begin{equation*}
-(N-p)\leq \frac{(h^2-(N-p)|\nabla h|^2)}{\rho^2} \leq 1.
\end{equation*}
Consequently,
\begin{equation}\label{eq:Qminus-evolution}
 \partial_\tau Q_{-}\leq q-1-N+C_{N,p,q}z -(N-p+1)z\left(h\sum_i b^{ii}
 +\frac{N-p}{h}\sum_i b^{ii}h_i^2\right).
\end{equation}
By $C^0$ and $C^1$ estimates, $z=\eta h^{p-1-N}e^{- Q_{-}}\leq C_Ne^{- Q_{-}}$. With this choice of $N$, dropping the last non-positive term in \eqref{eq:Qminus-evolution}, at $(\xi_0,\tau_0)$, we obtain
\begin{equation*}
0\leq \partial_\tau Q_{-}\leq q-1-N+C_{N,p,q}C_Ne^{- Q_{-}},
\end{equation*}
which is negative for sufficiently large $Q_{-}$, a contradiction. Thus $Q_{-}$ is uniformly bounded above on $\cC_\theta\times[0,T_0]$, with a bound independent of $T_0<T$. Since $T_0<T$ is arbitrary, $\cK$ is uniformly bounded below throughout $[0,T)$.

\end{proof}

Now we establish the following estimates for the principal curvatures.

\begin{lemma}\label{s4:lem-C2*}
Let $pq\geq 0$ and $\theta \in (0, \frac{\pi}{2})$. Suppose $h(\cdot, \tau)$ is a positive, even, smooth strictly convex solution to the flow \eqref{eq:scalar-flow}. Then there exists a positive constant $C$, depending only on $n$, $p$, $q$, $\theta$, $\min_{\mathcal{C}_\theta} f$, $\|f\|_{C^2(\mathcal{C}_\theta)}$ and $\Sigma_0$, such that the principal curvatures $\kappa_1, \ldots, \kappa_n$ of $\Sigma_{\tau}$ satisfy
\begin{equation}\label{eq:C2}
\frac{1}{C} \leq \kappa_i \leq C.
\end{equation}
\end{lemma}

\begin{proof}
In view of the uniform upper bound for the Gauss curvature $\mathcal{K}$ of $\Sigma_{\tau}$ (cf. \autoref{prop:Gauss-bounds}), it remains to establish a uniform lower bound for its principal curvatures, or equivalently, a uniform upper bound for its principal radii of curvature. For this purpose, we employ the maximum principle with the following auxiliary function
\begin{equation*}
Q=\log\sigma_1+Ah+B|\nabla u|^2,
\end{equation*}
where $u=h/\ell$, $B>0$ and $A<0$ will be chosen in that order. Fix an arbitrary $T_0\in(0,T)$ and take a maximum on $\cC_\theta\times[0,T_0]$. If it occurs at time zero, the initial data bound it. Otherwise its time coordinate is positive, and at that point $Q_\tau\geq0$.
\begin{claim}
For each fixed $B>0$, one has $Q_\mu<0$ on $\partial\cC_\theta$ at positive times once $A<0$ is sufficiently negative, where $A$ can be chosen to be independent of $B$.
\end{claim}
If the claim is true, we know that the maximum of $Q$ cannot be attained on $\partial\mathcal{C}_{\theta}$. To prove the claim, fix any point $\xi_0\in \partial \mathcal{C}_\theta$ and then choose a local orthonormal frame $\{e_i\}_{i=1}^n$ such that $e_n = \mu$ and $\{b_{ij}\}$ is diagonal at $\xi_0$.  Differentiating $h=\ell u$ and using \eqref{eq:ell}--\eqref{eq:u-boundary} yields
\begin{equation}\label{eq:mixed-boundary-radius}
h_{\alpha\mu}=b_{\alpha\mu}=0, \qquad 1\leq\alpha\leq n-1.
\end{equation}
Since $u_\mu=0$ and $\partial_\mu|\nabla u|^2=-2\cot\theta|\nabla u|^2$,
\begin{equation*}
Q_\mu=\frac{(\sigma_1)_\mu}{\sigma_1} +A\cot\theta\,h-2B\cot\theta|\nabla u|^2.
\end{equation*}
Let $F=\log\det b$.  We find
\begin{equation*}
\begin{aligned} F_\mu &=b^{\mu\mu}\nabla_\mu b_{\mu\mu} +\sum_{\alpha<n}b^{\alpha\alpha} \nabla_\mu b_{\alpha\alpha}, \\ (\sigma_1)_\mu &=\nabla_\mu b_{\mu\mu} +\sum_{\alpha<n}\nabla_\mu b_{\alpha\alpha}. \end{aligned}
\end{equation*}
Codazzi, \eqref{eq:mixed-boundary-radius}, and $\nabla_{e_\alpha}\mu=\cot\theta\,e_\alpha$ give
\begin{equation}\label{eq:normal-tangential-radius}
\nabla_\mu b_{\alpha\alpha} =\cot\theta\,(b_{\mu\mu}-b_{\alpha\alpha}).
\end{equation}
Using \eqref{eq:normal-tangential-radius}, we obtain
\begin{equation}
\begin{aligned}
(\sigma_1)_\mu &=b_{\mu\mu}F_\mu-\cot\theta\,b_{\mu\mu}^2 \sum_{\beta<n}b^{\beta\beta} +2(n-1)\cot\theta\,b_{\mu\mu} -\cot\theta\sum_{\beta<n}b_{\beta\beta}. \label{eq:sigma-normal}
\end{aligned}
\end{equation}
By \eqref{eq:scalar-flow},
\begin{equation}\label{eq:log-flow-speed}
\log(h-h_\tau)=-F+\log\eta+\log f+(n+1-q)\log\rho+p\log h.
\end{equation}
Furthermore, the definition of $\rho$, together with \eqref{eq:mixed-boundary-radius}, yields
\begin{equation*}
\frac{\rho_\mu}{\rho} =\frac{b_{\mu\mu}h_\mu}{\rho^2} =\cot\theta\frac{hb_{\mu\mu}}{\rho^2}.
\end{equation*}
Since $(h-h_\tau)_\mu=\cot\theta\,(h-h_\tau)$ and $\eta$ is spatially constant, differentiating \eqref{eq:log-flow-speed} in the $\mu$-direction therefore gives
\begin{equation*}
F_\mu=(p-1)\cot\theta+(\log f)_\mu +(n+1-q)\cot\theta\frac{hb_{\mu\mu}}{\rho^2}.
\end{equation*}
Thus the $C^0$, $C^1$ estimates imply $|F_\mu|\leq C(1+b_{\mu\mu})$.  Inserting this in \eqref{eq:sigma-normal}, using $b_{\mu\mu}/\sigma_1\leq1$, and discarding the last non-positive term in \eqref{eq:sigma-normal}, for constants $C_0,C_1,C_2>0$, there is
\begin{equation*}
Q_\mu\leq C_0+\frac{b_{\mu\mu}^2}{\sigma_1} \left(C_1-C_2\sum_{\beta<n}b^{\beta\beta}\right) +A\cot\theta\,h.
\end{equation*}
If $\sum_{\beta<n}b^{\beta\beta}\geq C_1/C_2$, then $Q_\mu<0$ provided that $A\leq A_1<0$ is sufficiently negative. Otherwise, $\sum_{\beta<n} b^{\beta\beta} \leq C_1/C_2$. Combining this with $b_{\beta\beta}\geq (\sum_{\gamma<n} b^{\gamma\gamma})^{-1}$ for $\beta<n$ and \autoref{prop:Gauss-bounds}, it follows that
\begin{equation*}
b_{\mu\mu} =\frac{\det b}{\prod_{\beta<n}b_{\beta\beta}}\leq C.
\end{equation*}
Thus $Q_\mu$ is also strictly negative provided that $A\leq A_2<0$ is sufficiently negative. Then the claim follows directly for $A\leq\min\{A_1,A_2\}$.

Thus, provided that $A$ satisfies the condition above, the chosen maximum point $(\xi_0,\tau_0)$ with $\tau_0>0$ must satisfy $\xi_0\in\cC_\theta\setminus\partial\cC_\theta$. Choose a local orthonormal frame $\{e_i\}_{i=1}^n$ around $(\xi_0,\tau_0)$ such that $\{b_{ij}\}(\xi_0,\tau_0)$ is diagonal, and hence $\{h_{ij}\}(\xi_0,\tau_0)$ is also diagonal.

Let
\begin{equation*}
\cL=b^{ij}\nabla_{ij}, \qquad \mathcal P=\cL-\frac1{h-h_\tau}\partial_\tau.
\end{equation*}

At $(\xi_0,\tau_0)$, we have $\mathcal P Q\leq0$. Put
\begin{equation*}
T_{ij}=\ell_i u_j+\ell_j u_i+u\delta_{ij}.
\end{equation*}
Using  $(x-y)^2\geq x^2/2-y^2$, \eqref{eq:ell}, $C^0$ and $C^1$ estimates, we have
\begin{equation*}
\begin{aligned}
b^{ij}u_{ki}u_{kj}&=\sum_{i,k}\frac{(b_{ki}-T_{ki})^2}{\ell^2b_{ii}}\geq \sum_{i}\frac{(b_{ii}-T_{ii})^2}{\ell^2b_{ii}}\\
&\geq \sum_{i}\frac{1}{b_{ii}}\left(\frac{b^2_{ii}}{2\ell^2}-\frac{T^2_{ii}}{\ell^2} \right)\geq c_{0}\sigma_1-C\sum_i b^{ii}.
\end{aligned}
\end{equation*}
Differentiating \eqref{eq:b-u}, one finds
 \begin{equation}
\begin{split}\label{bijk}
b_{ijk}&=\ell u_{ijk}+\ell_k u_{ij}+\ell_{ik}u_j +\ell_i u_{jk} +\ell_{jk}u_{i}+\ell_j u_{ik}+u_{k}\delta_{ij}\\
&=\ell (u_{kij}-u_k \delta_{ij}+u_j\delta_{ki})+\ell_k u_{ij}+\ell_{ik}u_j +\ell_i u_{jk} +\ell_{jk}u_{i}+\ell_j u_{ik}+u_{k}\delta_{ij}\\
&=\ell u_{kij}+\ell_k u_{ij}+u_{j}\delta_{ik}+\ell_{i}u_{jk}+(1-\ell) u_i \delta_{jk}+\ell_j u_{ik}+(1-\ell)u_k \delta_{ij},
\end{split}
\end{equation}
where we used $\ell_{ij}=(1-\ell)\delta_{ij}$.

Due to $F_k=b^{ij}b_{ijk}$, applying \eqref{bijk}, we have
\begin{equation*}
\ell b^{ij}u_{kij} =F_k-\ell_kb^{ij}u_{ij}-2b^{ij}\ell_i u_{jk} +(\ell-1)u_k\sum_i b^{ii}-(2-\ell)b^{kj}u_j.
\end{equation*}
Consequently,
\begin{equation}
\begin{aligned}\label{LU}
\cL|\nabla u|^2 ={}&\frac2\ell\left( \langle\nabla F,\nabla u\rangle -\langle\nabla\ell,\nabla u\rangle b^{ij}u_{ij} -2b^{ij}\ell_i u_{jk}u_k\right) \\
&+\frac{2(\ell-1)}\ell|\nabla u|^2\sum_i b^{ii} -\frac{2(2-\ell)}\ell\sum_i b^{ii}u_i^2 +2b^{ij}u_{ki}u_{kj}.
\end{aligned}
\end{equation}
Since $b=\ell\nabla^2u+T$, we have
\begin{equation*}
b^{ij}u_{ij}=\ell^{-1}(n-b^{ij}T_{ij}), \qquad b^{ij}\ell_i u_{jk}u_k =\ell^{-1}\langle\nabla\ell,\nabla u\rangle -\ell^{-1}b^{ij}\ell_iT_{jk}u_k.
\end{equation*}
Combining this with \eqref{LU}, we have
\begin{equation}\label{eq:L-grad-u}
\cL|\nabla u|^2 \geq2c_0\sigma_1+\frac2\ell\langle\nabla F,\nabla u\rangle -C\left(1+\sum_i b^{ii}\right).
\end{equation}

On the other hand, using $u_\tau=h_\tau/\ell$, we deduce that
\begin{equation*}
\begin{aligned}
\frac1{h-h_\tau}\partial_\tau|\nabla u|^2
={}&\frac2{h-h_\tau}|\nabla u|^2
-\frac2\ell\langle\nabla u,\nabla\log(h-h_\tau)\rangle+\frac2\ell\langle\nabla u,\nabla\log\ell\rangle.
\end{aligned}
\end{equation*}
Now set
\begin{equation}\label{eq:psiG}
\begin{aligned}
\psi=F+\log(h-h_\tau)=\log\eta+p\log h+G,
\end{aligned}
\end{equation}
where
\begin{equation}\label{eq:psiG1}
\begin{aligned}
G=\log f+(n+1-q)\log\rho.
\end{aligned}
\end{equation}
From \eqref{eq:L-grad-u}, \eqref{eq:psiG} and \eqref{eq:psiG1}, one gets
\begin{equation}\label{eq:P-grad-u}
\mathcal P|\nabla u|^2 \geq2c_0\sigma_1+\frac2\ell\langle\nabla\psi,\nabla u\rangle -C\left(1+\sum_i b^{ii}\right).
\end{equation}

Since $b=\nabla^2h+hI$ is a Codazzi tensor on the unit sphere, the contracted commutation identity is
\begin{equation*}
\Delta b_{ij}=\nabla_{ij}^2\sigma_1+nb_{ij}-\sigma_1\delta_{ij}.
\end{equation*}
Thus,
\begin{equation}\label{SIE}
\mathcal{L}\sigma_1=b^{ij}\Delta b_{ij}-n^2+\sigma_1\sum_i b^{ii}.
\end{equation}
Twice differentiating $F=\log\det b$ yields
\begin{equation}\label{eq:Delta-F-full}
\Delta F=b^{ij}\Delta b_{ij} -b^{ir}b^{sj}\nabla_kb_{rs}\nabla_kb_{ij}.
\end{equation}
Substituting \eqref{eq:Delta-F-full} into \eqref{SIE}, at $(\xi_0, \tau_0)$, we obtain
\begin{equation*}
\mathcal{L}\sigma_1=\Delta F+\sigma_1 \sum_i b^{ii}-n^2+\sum_{i,j,k}b^{ii}b^{kk}(\nabla_{j}b_{ik})^2.
\end{equation*}
It follows that
\begin{equation}\label{SIE2}
\mathcal{L} \log\sigma_1=\frac{\mathcal{L} \sigma_1}{\sigma_1}-\frac{b^{ij}(\sigma_1)_i (\sigma_1)_j}{\sigma^2_1}\geq \frac{\Delta F}{\sigma_1}+\sum_i b^{ii}-\frac{n^2}{\sigma_1},
\end{equation}
where we used that for each $i$,
\begin{equation*}
\sum_{j,k}b^{kk}(\nabla_jb_{ik})^2
\geq\sum_jb^{jj}(\nabla_i b_{jj})^2
\geq\frac{(\nabla_i\sigma_1)^2}{\sigma_1}.
\end{equation*}
By a direct calculation, we have
\begin{equation*}
(\sigma_1)_\tau=\sigma_1-\bigl(\Delta(h-h_\tau)+n(h-h_\tau)\bigr),
\end{equation*}
and
\begin{equation*}
\begin{aligned}
-\frac1{h-h_\tau}\partial_\tau\log\sigma_1
&=-\frac1{h-h_\tau}+\frac{\Delta(h-h_\tau)}{(h-h_\tau)\sigma_1}+\frac n{\sigma_1},\\
\frac{\Delta(h-h_\tau)}{h-h_\tau}
&=\Delta\log(h-h_\tau)+|\nabla\log(h-h_\tau)|^2.
\end{aligned}
\end{equation*}
Combining this with \eqref{SIE2} and \eqref{eq:psiG}, we obtain
\begin{equation}\label{eq:P-logsigma}
\begin{aligned}
\mathcal P\log\sigma_1\geq{}&\frac{\Delta\psi}{\sigma_1}+\sum_i b^{ii}
-\frac{n(n-1)}{\sigma_1}-\frac1{h-h_\tau}+\frac{|\nabla\log(h-h_\tau)|^2}{\sigma_1}.
\end{aligned}
\end{equation}

Finally, since $\cL h=n-h\sum_i b^{ii}$ and $h_\tau/(h-h_\tau)=h/(h-h_\tau)-1$,
\begin{equation}\label{eq:P-h}
\mathcal P h=n+1-h\sum_i b^{ii}-\frac h{h-h_\tau}.
\end{equation}

Since $\eta$ is spatially constant,
\begin{equation*}
\Delta\psi=p\Delta\log h+\Delta G,
\end{equation*}
and
\begin{equation*}
\frac{\Delta\log h}{\sigma_1} =\frac1h-\frac n{\sigma_1} -\frac{|\nabla h|^2}{h^2\sigma_1}.
\end{equation*}
By \autoref{prop:Gauss-bounds}, $C^0$ and $C^1$ estimates, one gets
\begin{equation*}
\frac{1}{C}\leq h-h_\tau=\eta f \rho^{n+1-q}h^p\cK\leq C.
\end{equation*}

Combining these identities with \eqref{eq:P-logsigma}, \eqref{eq:P-h}, and \eqref{eq:P-grad-u} gives
\begin{equation}\label{eq:PQ-first}
\mathcal P Q \geq2c_0B\sigma_1+(-Ah-CB)\sum_i b^{ii} -C\left(1+|A|+B+\frac1{\sigma_1}\right)+\Psi,
\end{equation}
where
\begin{equation}\label{eq:Psi}
\Psi=\frac{\Delta G}{\sigma_1} +\frac{2B}{\ell}\langle\nabla u,\nabla G\rangle.
\end{equation}

At $(\xi_0,\tau_0)$, we have
\begin{equation}\label{eq:gradQzero}
\frac{\nabla\sigma_1}{\sigma_1} =-A\nabla h-B\nabla|\nabla u|^2.
\end{equation}
Using $\rho^2=h^2+|\nabla h|^2$, we have $\nabla(\rho^2)=2b\nabla h$ and
\begin{equation}
\begin{aligned}
\Delta(\rho^2) &=2h\sigma_1-2nh^2+2|\nabla^2h|^2 +2\langle\nabla h,\nabla\sigma_1\rangle.\label{eq:Delta-rho2}
\end{aligned}
\end{equation}
Since $G=\log f+\frac{n+1-q}{2}\log(\rho^2)$,
\begin{equation*}
\Delta G=\Delta\log f+\frac{n+1-q}{2} \left(\frac{\Delta(\rho^2)}{\rho^2} -\frac{|\nabla(\rho^2)|^2}{\rho^4}\right).
\end{equation*}
Substituting \eqref{eq:gradQzero} and \eqref{eq:Delta-rho2} into \eqref{eq:Psi} gives
\begin{equation*}
\begin{aligned}
\Psi={}&\frac{\Delta\log f}{\sigma_ 1} +\frac{2B}{\ell}\langle\nabla u,\nabla\log f\rangle \\
&+(n+1-q)\left( \frac h{\rho^2}-\frac{nh^2}{\sigma_1\rho^2} +\frac{|\nabla^2h|^2}{\sigma_1\rho^2} -A\frac{|\nabla h|^2}{\rho^2}\right) \\
&-\frac{n+1-q}{2}\frac{|\nabla(\rho^2)|^2}{\sigma_1\rho^4} +\frac{(n+1-q)B}{\rho^2}\chi,
\end{aligned}
\end{equation*}
where
\begin{equation*}
\chi=\ell^{-1}\langle\nabla u,\nabla(\rho^2)\rangle -\langle\nabla h,\nabla|\nabla u|^2\rangle.
\end{equation*}
The potentially unbounded Hessian terms in $\chi$ cancel. Using $\nabla(\rho^2)=2b\nabla h$ and $b=\ell\nabla^2u+T$, one has
\begin{equation*}
\begin{aligned}
\ell^{-1}\langle\nabla u,\nabla(\rho^2)\rangle &=\frac2\ell\langle\nabla u,b\nabla h\rangle, \\
\langle\nabla h,\nabla|\nabla u|^2\rangle &=2\langle\nabla h,\nabla^2u\,\nabla u\rangle =\frac2\ell\langle\nabla h,(b-T)\nabla u\rangle.
\end{aligned}
\end{equation*}
Therefore
\begin{equation*}
\chi=\frac2\ell\langle\nabla h,T\nabla u\rangle,
\end{equation*}
which is uniformly bounded by the $C^0$ and $C^1$ estimates.  Furthermore,
\begin{equation*}
|\nabla(\rho^2)|^2\leq C(1+|\nabla^2h|^2), \qquad \frac{|\nabla^2h|^2}{\sigma_1} \leq\sigma_1+\frac{nh^2}{\sigma_1}.
\end{equation*}
From the above, we further have
\begin{equation}
\begin{split}
\label{eq:Psi-lower}
\Psi\geq-C_0\sigma_1 -C_1\left(1+|A|+B+\frac1{\sigma_1}\right).
   \end{split}
\end{equation}

Choose $B>0$ large enough so that $2c_0B>C_0+1$. Let $c_h>0$ be the time-independent lower bound for $h$ in \eqref{eq:C0}. Choose $A<0$ sufficiently negative that
\begin{equation*}
-Ac_h\geq CB+1,\qquad A\leq\min\{A_1,A_2\}<0.
\end{equation*}
These choices are independent of $T_0<T$. Then at $(\xi_0,\tau_0)$, \eqref{eq:PQ-first} and \eqref{eq:Psi-lower} imply
\begin{equation*}
\begin{split}
0\geq \mathcal P Q &\geq \sigma_1-C\left(1+\frac1{\sigma_1}\right)\\
&= e^{Q-Ah-B|\nabla u|^2}-\frac{C}{e^{Q-Ah-B|\nabla u|^2}} -C.
\end{split}
\end{equation*}
Thus $Q(\xi_0,\tau_0)\leq C$. The boundedness of $Ah+B|\nabla u|^2$ and maximality of $Q(\xi_0,\tau_0)$ imply $\sigma_1\leq C$ on $\cC_\theta\times[0,T_0]$, with $C$ independent of $T_0<T$. Since $T_0<T$ is arbitrary, this bound holds throughout $\cC_\theta\times[0,T)$. Thus $b\leq CI$. The lower bound $\det b=\cK^{-1}\geq c$ then gives $b\geq c\,C^{1-n}I$, proving the two-sided principal curvature bounds.
\end{proof}

\section{Long-Time Existence and Convergence}\label{Sec5}

On the boundary $\partial \cC_\theta$, write $\mathcal Bv:=\nabla_\mu v-\cot\theta\,v$. We first establish short-time existence. 
\begin{prop}\label{prop:short-time}
Let $h_0\in C^\infty(\cC_\theta)$ be positive and strictly convex, with $\mathcal Bh_0=0$. Then, for some $T_{\mathrm{loc}}>0$ and $\alpha\in(0,1)$, \eqref{eq:scalar-flow} has a unique, positive strictly convex solution $h(\cdot, \tau)$ in $C^{2+\alpha,1+\alpha/2}(\cC_\theta\times[0,T_{\mathrm{loc}}])$.  
Moreover, if $f$ and $h_0$ are even, then $h(\cdot,\tau)$ remains even  as long as the solution exists.
\end{prop}

\begin{proof}

Set $\mathcal F[\tilde h]=f(\tilde  h^2+|\nabla\tilde  h|^2)^{(n+1-q)/2}\tilde  h^p(\det b[\tilde  h])^{-1}$. Here $b[v]=\nabla^2v+vI$ and $D\mathcal F$ denotes the Fr\'echet derivative. We first solve $\tilde h_t=-\mathcal F[\tilde  h]$ with $\mathcal B\tilde  h=0$ and initial value $h_0$. The linearized operator $L_0=-D\mathcal F[h_0]$ has a uniformly positive definite principal coefficient matrix $\mathcal F[h_0](b[h_0]^{ij})$. Consequently, combining the linear estimates for oblique boundary value problems in \cite[Theorem~5.18]{Lie96} with a standard fixed-point argument based on the contraction mapping theorem, we obtain a unique local solution in $C^{2+\alpha,1+\alpha/2}$.

Since $\eta[\tilde  h]>0$, define 
\begin{equation*}
\tau(t)=\int_0^t\frac{ds}{\eta[\tilde  h(\cdot,s)]},\qquad
h(\xi,\tau(t))=e^{\tau(t)}\tilde  h(\xi,t).
\end{equation*}
Using $\mathcal F[\lambda\tilde  h]=\lambda^{p+1-q}\mathcal F[\tilde  h]$ and $\eta[\lambda\tilde  h]=\lambda^{q-p}\eta[\tilde  h]$, we have $h_\tau=h-\eta[h]\mathcal F[h]$. Conversely, the transformation
\begin{equation*}
t(\tau)=\int_0^\tau e^{(p-q)s}\eta[h(\cdot,s)]\,ds,\qquad
\tilde  h(\xi,t(\tau))=e^{-\tau}h(\xi,\tau)
\end{equation*}
recovers the unnormalized equation. The existence and uniqueness of local solutions to \eqref{eq:scalar-flow} follow from those for the unnormalized equation.  For the horizontal reflection $\mathcal{R}(\xi',\xi_{n+1})=(-\xi',\xi_{n+1})$, evenness of $f,h_0$ implies that $h(\mathcal{R}\xi,\tau)$ solves the same equation, so evenness is preserved by uniqueness.
\end{proof}

We next use parabolic boundary regularity to obtain higher-order
bounds from the previous a priori estimates, and then apply them to prove long-time existence. Let $T_\ast$ be the maximal existence time of the positive strictly convex solution to \eqref{eq:scalar-flow}. By \eqref{eq:C0}, \eqref{eq:C1}, \eqref{Bound:K} and \eqref{eq:C2}, we have
\begin{equation*}
cI\leq b\leq CI,\qquad c\leq h-h_\tau\leq C,\qquad |h_\tau|\leq C.
\end{equation*}
This implies that \eqref{eq:scalar-flow} is uniformly parabolic and concave in the Hessian variable. Moreover, by \autoref{prop:preserve}, we have
\begin{equation*}
\left|\frac{\eta_\tau}{\eta}\right|
=|p|\left|\frac{\int_{\cC_\theta}fh^{p-1}h_\tau\,d\xi}{H_p}\right|\leq C.
\end{equation*}
Thus $|\eta_\tau|\leq C$, and spatial interpolation gives uniform Lipschitz continuity of $h$ and $1/2$-H\"older continuity of $\nabla h$ in time. Now let $w=h_t$. Differentiating \eqref{eq:scalar-flow} with respect to time, we find that $w$ satisfies a linear, uniformly parabolic equation whose coefficients, including those of the second-order, first-order, and zeroth-order terms, are uniformly bounded by the previously established estimates. Parabolic H\"older estimates for Neumann boundary conditions  then yield a uniform H\"older bound for $w$. Combining this bound with
a bootstrap argument based on the regularity theory for parabolic
equations with Neumann boundary conditions
(see, e.g., \cite[Theorems~6.1, 6.4, and~6.5]{D88}
and \cite[Theorem~14.23]{Lie96}), for $0<\tau_0<T_\ast$, we obtain
\begin{equation}\label{eq:higher}
\sup_{\tau_0\leq\tau<T_\ast}
\left(\|h(\cdot,\tau)\|_{C^k(\cC_\theta)}+\|h_\tau(\cdot,\tau)\|_{C^k(\cC_\theta)}\right)
\leq C(k,\tau_0),\qquad {\rm for\ every \ integer} \  k\geq 0.
\end{equation}
If $T_\ast<\infty$, \eqref{eq:higher} implies that $h(\cdot,\tau)$ converges to a smooth limit $h_*$ in  every $C^k$ norm as $\tau\rightarrow T_\ast$. Moreover, $h_*$ satisfies $h_*\geq c$, $b[h_*]\geq cI$ and $\mathcal Bh_*=0$. Proposition~\ref{prop:short-time} therefore extends the solution
beyond $T_\ast$, contradicting maximality. Hence
$T_\ast=\infty$.

We now recover the flow \eqref{eq:geometric-flow} from the scalar version \eqref{eq:scalar-flow}
by a tangential reparametrization that preserves the boundary. Set $Y=\nabla h+h\nu$, where $\nu(\xi)=\xi-\cos\theta\,e$, and define
\begin{equation*}
a^i=\ell b^{ij}\nabla_j\left(\frac{h-h_\tau}{\ell}\right).
\end{equation*}
Using $dY(a)=b_{ij}a^j e_i$ and $\widetilde\nu=\ell\nu+\nabla\ell$, we have
\begin{equation*}
Y_\tau+dY(a)=Y-\frac{h-h_\tau}{\ell}\widetilde\nu.
\end{equation*}
Since $\mathcal Bh=\mathcal Bh_\tau=\mathcal B\ell=0$ and $b_{\alpha\mu}=0$,
\begin{equation*}
\nabla_\mu\left(\frac{h-h_\tau}{\ell}\right)=0,
\qquad a^\mu=0
\quad\text{on }\partial\cC_\theta.
\end{equation*}
Let $\varphi$
solve $\varphi_\tau=a\circ\varphi$ for $\tau>0$, then  $X=Y\circ\varphi$  satisfies
$\partial_\tau X
= X-\left(\frac{h-h_\tau}{\ell}\right)\circ\varphi\,
\widetilde{\nu}\circ\varphi$. By applying \eqref{eq:scalar-flow}, we obtain \eqref{eq:geometric-flow}.

Using the above results, we now prove  \autoref{thm:main}. The uniform a priori estimates established in \autoref{Sec4} imply that $\Phi_{p,q}(\Sigma_\tau)$ remains bounded for $\tau\geq0$. Together with \autoref{prop:entropy}, this shows
\begin{equation*}
\int_0^\infty\|h_\tau(\cdot,\tau)\|_{L^2(\cC_\theta)}^2\,d\tau<\infty,
\end{equation*}
 which implies that there exists a sequence $\tau_j\to\infty$ such that $||h_{\tau}(\cdot,\tau_j)||_{L^2(\mathcal{C}(\theta))}\rightarrow 0$. By \eqref{eq:higher} and the Arzel\`a-Ascoli theorem, we may pass to a subsequence, still denoted by $\{\tau_j\}$, such that $h(\cdot,\tau_j)$ converges smoothly to a positive, strictly convex function $h_\infty$, which remains even and satisfies the Robin boundary condition. Continuity of $Q_q$ and $H_p$ implies $\eta[h(\cdot,\tau_j)]\to\eta[h_\infty]$. Consequently, $h_{\infty}$ satisfies \eqref{eq:normalized-target}. Moreover, the capillary hypersurface $\Sigma$ determined by the capillary support function $h_{\infty}$ is smooth, even, and strictly convex. The two sides of \eqref{eq:target} are homogeneous of degrees $n$ and $n+p-q$, respectively. If $p\neq q$, the dilation $\bar h=\eta[h_\infty]^{1/(p-q)}h_\infty$ solves \eqref{eq:target}. This completes the proof of \autoref{thm:main}.

\end{document}